\documentclass[12pt]{article}

\begin{document}

\begin{titlepage}
\title{\bf Constrained Paracomplex Mechanical Systems}
\author{ Mehmet Tekkoyun \footnote{tekkoyun@pau.edu.tr} \\
 {\small Department of Mathematics, Pamukkale University,}\\
{\small 20070 Denizli, Turkey}\\
Murat Sari \footnote{msari@pau.edu.tr} \\
 {\small Department of Mathematics, Pamukkale University,}\\
{\small 20070 Denizli, Turkey}}
\date{\today}
\maketitle

\begin{abstract}
In this study, it is introduced paracomplex analogue of
Lagrangians and Hamiltonians with constraints in the framework of
para-K\"{a}hlerian manifolds. The geometrical and mechanical
results on the constrained mechanical system have also been discussed.


\end{abstract}

\section{Introduction}

It is well known that the dynamics of Lagrangian and Hamiltonian formalisms
is characterized by a suitable vector field defined on the tangent and
cotangent bundles which are phase-spaces of velocities and momentum of a
given configuration manifold. If $Q$ is an $m$-dimensional configuration
manifold and $L:TQ\rightarrow \mathbf{R}$ is a regular Lagrangian function,
then there is a unique vector field $\xi _{L}$ on $TQ$ such that dynamical
equations
\begin{equation}
i_{\xi _{L}}\Phi _{L}=dE_{L},  \label{1.1}
\end{equation}
where $\Phi _{L}$ is the symplectic form and $E_{L}$ is the energy
associated to $L$. The Euler-Lagrange vector field $\xi _{L}$ is a second
order differential equation on $Q$ since its integral curves are the
solutions of the Euler-Lagrange equations given by
\begin{equation}
\frac{d}{dt}\frac{\partial L}{\partial \dot{q}^{i}}-\frac{\partial L}{%
\partial q^{i}}=0,  \label{1.2}
\end{equation}
where $q^{i}$ and $(q^{i},{\dot{q}}^{i}),1\leq i\leq m,$ are coordinate
systems of $Q$ and $TQ,$ respectively$.$ The triple $(TQ,\Phi _{L},E_{L})$
is called \textit{Lagrangian mechanical system} on the tangent bundle $TQ.$
Assume that $(TQ,\Phi _{L})$ is symplectic manifold and $\overline{\omega }%
=\{\omega _{1},...,\omega _{r}\}$ is a system of constraints on $TQ.$ It is
called as a \textit{constraint} on $TQ$\ to a non-zero 1-form $\omega
=\wedge ^{a}\omega _{a}$ on $TQ,$ such that $\wedge ^{a}$, $1\leq a\leq r,$
are Lagrange multipliers. The quartet $(TQ,\Phi _{L},E_{L},\overline{\omega }%
)$ is said to be a \textit{regular Lagrangian system with constraints}. The
constraints $\overline{\omega }$ are said to be classical constraints if the
1-forms $\omega _{a}$ are basic. Then holonomic classical constraints define
foliations on the configuration manifold $Q,$ but holonomic constraints also
admit foliations on the phase space of velocities $TQ.$ As is the case in
real studies, generally, a curve $\alpha $ satisfying the Euler Lagrange
equations for Lagrangian energy $E_{L}$ does not satisfy the constraints. In
order to satisfy the constraints, some additional forces act on the system
as well as \textit{force} $dE_{L}$ for a curve $\alpha $ $.$ It is said that
the quartet $(TQ,\Phi _{L},E_{L},\overline{\omega })$ defines a \textit{%
mechanical system with constraints} if vector field $\xi $ given by the
equations of motion
\begin{equation}
i_{\xi }\Phi _{L}=dE_{L}+\wedge ^{a}\omega _{a},\,\,\,\,\,\,\omega _{a}(\xi
)=0,  \label{1.3}
\end{equation}
is a second order differential equation. Then, it is given Euler-Lagrange
equations with constraints as follows:
\begin{equation}
\frac{\partial L}{\partial q^{i}}-\frac{d}{dt}\frac{\partial L}{\partial
\stackrel{.}{q}^{i}}=\wedge ^{a}(\omega _{a})_{i},  \label{1.4}
\end{equation}

If $H:T^{*}Q\rightarrow \mathbf{R}$ is a regular Hamiltonian function then
there is a unique vector field $Z_{H}$ on cotangent bundle $T^{*}Q$ such
that dynamical equations
\begin{equation}
i_{Z_{H}}\Phi =dH,  \label{1.5}
\end{equation}
where $\Phi $ is the symplectic form and $H$ stands for Hamiltonian
function. The paths of the Hamiltonian vector field $Z_{H}$ are the
solutions of the Hamiltonian equations shown by

\begin{equation}
\frac{dq^{i}}{dt}=\frac{\partial H}{\partial p_{i}},\,\,\,\,\,\,\,\,\,\frac{%
dp_{i}}{dt}=-\frac{\partial H}{\partial q^{i}},  \label{1.6}
\end{equation}
where $q^{i}$ and $(q^{i},p_{i}),1\leq i\leq m,$ are coordinates of $Q$ and $%
T^{*}Q.$ The triple $(T^{*}Q,\Phi ,H),$ is called \textit{Hamiltonian system}
on the cotangent bundle $T^{*}Q$ with symplectic form $\Phi $. Let $T^{*}Q$%
\textbf{\ }be symplectic manifold with\textbf{\ }closed symplectic form $%
\Phi .$ Similar to constraints on $TQ,$ by a \textit{constraint} on $T^{*}Q$
is said to be a \textit{non-zero 1-form} $\omega =\wedge ^{a}\omega _{a}$ on
$T^{*}Q.$ A set $\overline{\omega }=\{\omega _{1},...,\omega _{s}\}$ of $s$
linearly independent 1- forms on $T^{*}Q$ may be named to be a \textit{%
system of constraints} on $T^{*}Q$. We say that a curve $\alpha $ in $T^{*}Q$
satisfies the constraints if $\omega _{a}(\stackrel{.}{\alpha }%
(t))=0,\,1\leq a\leq s.$

Let $(T^{*}Q,\Phi ,H)$ be a Hamiltonian system\textit{\ }on\textbf{\ }%
symplectic manifold $T^{*}Q$ with\textbf{\ }closed symplectic form $\Phi $.
Let us consider a Hamiltonian system $(T^{*}Q,\Phi ,H)$ together with a
system $\overline{\omega }$ of constraints on $T^{*}Q.$ So, it is called $%
(T^{*}Q,\Phi ,H,\overline{\omega })$ to be a \textit{Hamiltonian system with}
\textit{constraints. }In\textit{\ }general, a curve $\alpha $ satisfying the
Hamiltonian equations for energy $H$ does not satisfy the constraints. For a
curve $\alpha $ satisfying the constraints, some additional forces must act
on the system in addition to the \textit{force} $dH.$ So, the dynamical
equations of motion become
\begin{equation}
i_{Z}\Phi =dH+\wedge ^{a}\omega _{a},\,\,\,\,\,\,\,\,\,\omega _{a}(Z)=0,
\label{1.7}
\end{equation}
where $Z$ is a vector field on $T^{*}Q$. From (\ref{1.7}), Hamiltonian
equations with constraints is given by:

\begin{equation}
\begin{array}{l}
\frac{dq^{i}}{dt}=(\frac{\partial H}{\partial p_{i}}+\wedge ^{a}(B_{a})_{i}),
\\
\frac{dp_{i}}{dt}=\mathbf{-}(\frac{\partial H}{\partial q_{i}}+\wedge
^{a}(A_{a})_{i}), \\
(A_{a})_{i}\frac{dq^{i}}{dt}+(B_{a})_{i}\frac{dp_{i}}{dt}=0,
\end{array}
\label{1.8}
\end{equation}
where $1\leq i\leq m,\,1\leq a\leq s.$

It can be easily understood that the above approach provides a good
framework for studying Lagrangian and Hamiltonian formalisms of classical
mechanics. There are some articles in \cite{mcrampin, nutku,weber,tekkoyun}
and books in \cite{deleon1985,deleon} on differential geometric methods in
mechanics. It is well known that (para)K\"{a}hlerian manifolds play an
essential role in various areas of mathematics and mathematical physics, in
particular, in the theory of dynamical systems, algebraic geometry, the
geometry of Einstein manifolds, quantum mechanics, quantum field theory, and
in the theory of superstrings and nonlinear sigma-models, too. For example,
it was shown in \cite{Ivanov} that the reflector space of an Einstein
self-dual non-Ricci flat 4-manifold as well as the reflector space of a
paraquaternionic K\"{a}hlerian manifold admit both Nearly para-K\"{a}hlerian
and almost para-K\"{a}hlerian structures. Wade \cite{Wade} showed that
generalized paracomplex structures are in one-to-one correspondence with
pairs of transversal Dirac structures on a smooth manifold. In \cite{Ulrych}%
, it was given a representation of the quadratic Dirac equation and the
Maxwell equations in terms of the three-dimensional universal complex
Clifford algebra C$_{3,0}$. Baylis and Jones introduced in \cite{baylis1}
that a R$_{3,0}$ Clifford algebra has enough structure to describe
relativity as well as the more usual R$_{1,3}$ Dirac algebra or the R$_{3,1}$
Majorana algebra. In \cite{baylis2}, Baylis represented relativistic
space-time points as paravectors and applies these paravectors to
electrodynamics. Tekkoyun \cite{tekkoyun1} generalized the concept of
Hamiltonian dynamics with constraints to complex case. In the above studies;
although paracomplex geometry, complex mechanical systems with constraints,
Lagrangian and Hamiltonian mechanics were given in a tidy and nice way, they
have not dealt with constrained paracomplex mechanical systems. Therefore,
in this paper, as a contribution to the modern development of Lagrangian and
Hamiltonian systems of classical mechanics, it was obtained paracomplex
analogous of some topics in the geometric theory of constraints given in
\cite{weber,deleon,tekkoyun1}, and it has an important role in mechanical
systems as pointed out in the above.

The present paper is structured as follows. In sections 1, 2 and
3, it is recalled paracomplex, para-Hermitian and
para-K\"{a}hlerian manifolds, and also para-Euler-Lagrange
equations and para-Hamiltonian equations on para-K\"{a}hlerian
manifolds. In sections 4 and 5, paracomplex Euler-Lagrange and
Hamiltonian equations with constraints on para-K\"{a}hlerian
manifold are deduced. In the conclusion section, the geometrical
and mechanical theory of para mechanical system with constraints
was presented.

\section{Preliminaries}

In this paper, all the geometrical objects are differentiable and the
Einstein summation convention is in use. So, $\mathbf{A}$, $\mathcal{F}(TM)$%
, $\chi (TM)$ and $\Lambda ^{1}(TM)$ denote the set of paracomplex numbers,
the set of paracomplex functions on tangent bundle $TM$, the set of
paracomplex vector fields on tangent bundle $TM$ and the set of paracomplex
1-forms on tangent bundle $TM$, respectively. Here $1\leq i\leq m.$ Some
geometric structures on the differential manifold $M$ given by \cite
{cruceanu} can be extended to $TM$ as follows:

\subsection{Paracomplex Geometry}

An \textit{almost product structure} $J$ on a tangent bundle $TM$ of $m$%
-real dimensional configuration manifold $M$ is a (1,1) tensor field $J$ on $%
TM$ such that $J^{2}=I.$ Here, the pair $(TM,J)$ is called an \textit{almost
product manifold}. An \textit{almost paracomplex manifold} is an almost
product manifold $(TM,J)$ such that the two eigenbundles $TT^{+}M$ and $%
TT^{-}M$ associated to the eigenvalues +1 and -1 of $J$, respectively, have
the same rank. The dimension of an almost paracomplex manifold is
necessarily even. Equivalently, a splitting of the tangent bundle $TTM$ of
tangent bundle $TM$, into the Whitney sum of two subbundles on $TT^{\pm }M$
of the same fiber dimension is called an \textit{almost paracomplex structure%
} on $TM.$ From physical point of view, this splitting means that a
reference frame has been chosen. Obviously, such a splitting is broken under
reference frame transformations. An almost paracomplex structure on a $2m$%
-dimensional manifold $TM$ may alternatively be defined as a $G$-structure
on $TM$ with structural group $GL(n,\mathbf{R})\times GL(n,\mathbf{R})$.

A \textit{paracomplex manifold is }an almost paracomplex manifold $(TM,J)$
such that $G$- structure defined by tensor field $J$ is integrable. Let $%
(x^{i})$ and $(x^{i},\,y^{i})$ be a real coordinate system of $M$ and $TM,$
and $\{(\frac{\partial }{\partial x^{i}})_{p},(\frac{\partial }{\partial
y^{i}})_{p}\}$ and $\{(dx^{i})_{p},(dy^{i})_{p}\}$ natural bases over $%
\mathbf{R}$ of tangent space $T_{p}(TM)$ and cotangent space $T_{p}^{*}(TM)$
of $TM,$ respectively$.$ Then, $J$ can be denoted as
\begin{equation}
J(\frac{\partial }{\partial x^{i}})=\frac{\partial }{\partial y^{i}},\,J(%
\frac{\partial }{\partial y^{i}})=\frac{\partial }{\partial x^{i}}.
\label{2.1}
\end{equation}
Let $\ z^{i}=\ x^{i}+$\textbf{j}$\ y^{i},$ \textbf{j}$^{2}=1,$ be a
paracomplex local coordinate system of $TM.$ The vector and covector fields
are defined, respectively, as follows:
\begin{equation}
(\frac{\partial }{\partial z^{i}})_{p}=\frac{1}{2}\{(\frac{\partial }{%
\partial x^{i}})_{p}-\mathbf{j}(\frac{\partial }{\partial y^{i}})_{p}\},\,(%
\frac{\partial }{\partial \overline{z}^{i}})_{p}=\frac{1}{2}\{(\frac{%
\partial }{\partial x^{i}})_{p}+\mathbf{j}(\frac{\partial }{\partial y^{i}}%
)_{p}\},  \label{2.2}
\end{equation}
\begin{equation}
\left( dz^{i}\right) _{p}=\left( dx^{i}\right) _{p}+\mathbf{j}%
(dy^{i})_{p},\,\left( d\overline{z}^{i}\right) _{p}=\left( dx^{i}\right)
_{p}-\mathbf{j}(dy^{i})_{p}.  \label{2.3}
\end{equation}
The above equations represent the bases of tangent space $T_{p}(TM)$ and
cotangent space $T_{p}^{*}(TM)$ of $TM$, respectively. Then the following
results can be easily obtained, respectively:
\begin{equation}
J(\frac{\partial }{\partial z^{i}})=-\mathbf{j}\frac{\partial }{\partial
z^{i}},\,J(\frac{\partial }{\partial \overline{z}^{i}})=\mathbf{j}\frac{%
\partial }{\partial \overline{z}^{i}},  \label{2.4}
\end{equation}
\begin{equation}
J^{*}(dz^{i})=-\mathbf{j}dz^{i},\,J^{*}(d\overline{z}^{i})=\mathbf{j}d%
\overline{z}^{i}.  \label{2.5}
\end{equation}
Here, $J^{*}$ stands for the dual endomorphism of cotangent space $%
T_{p}^{*}(TM)$ of manifold $TM$ satisfying $J^{*2}=I$ .

An \textit{almost para}-\textit{Hermitian manifold }$(TM,g,J)$ is a
differentiable manifold $TM$ endowed with an almost product structure $J$
and a pseudo-Riemannian metric $g$, compatible in the sense that
\begin{equation}
g(JX,Y)+g(X,JY)=0,\,\,\,\,\forall X,Y\in \chi (TM).  \label{2.6}
\end{equation}

An \textit{almost para}-\textit{Hermitian structure }on a differentiable
manifold\textit{\ }$TM$\textit{\ }is $G$-structure on $TM$ whose structural
group is the representation of the paraunitary group $U(n,\mathbf{A})$ given
in \cite{cruceanu}. A \textit{para}-\textit{Hermitian manifold }is a
manifold with an integrable almost para-Hermitian structure $(g,J)$.
2-covariant skew tensor field $\Phi $ defined by $\Phi (X,Y)=g(X,JY)$ is
so-called as \textit{fundamental 2-form}$.$ An almost para-Hermitian
manifold $(TM,g,J)$, such that $\Phi $ is closed, is so-called as an \textit{%
almost para-K\"{a}hlerian manifold}.

A para-Hermitian manifold $(TM,g,J)$ is said to be a \textit{%
para-K\"{a}hlerian manifold} if $\Phi $ is closed. Also, by means of
geometric structures, one may show that $(T^{*}M,g,J)$ is a \textit{%
para-K\"{a}hlerian manifold.}

\subsection{Paracomplex Lagrangian Systems}

In this section, some paracomplex fundamental concepts and
para-Euler-Lagrange equations for classical mechanics structured on
para-K\"{a}hlerian manifold $TM$ introduced in \cite{tekkoyun} can be
recalled.

Let $J$ be an almost paracomplex structure on the para-K\"{a}hlerian
manifold and $(z^{i},\overline{z}^{i})$ its coordinates. Let a second order
differential equation be vector field $\xi _{L}$ given by:
\begin{equation}
\xi _{L}=\xi ^{i}\frac{\partial }{\partial z^{i}}+\overline{\xi }^{i}\frac{%
\partial }{\partial \overline{z}^{i}},  \label{2.7}
\end{equation}
Then vector field $V=J\xi _{L}$ is called a \textit{para}-\textit{Liouville
vector field} on the para-K\"{a}hlerian manifold $TM$. The mappings given by
$T,P:TM\rightarrow \mathbf{A}$ such that $T=\frac{1}{2}m_{i}(\stackrel{.}{z}%
^{i})^{2},\,\,P=m_{i}\mathbf{g}h$ can be called as \textit{the kinetic energy%
} and \textit{the potential energy of system,} respectively,\textit{\ }where%
\textit{\ }$m_{i}$ is mass of a mechanical system, $\mathbf{g}$ is the
gravity and $h$ is the distance of the mechanical system on the
para-K\"{a}hlerian manifold to the origin. Then we call map $L:TM\rightarrow
\mathbf{A}$ such that $L=T-P$ as \textit{para}-\textit{Lagrangian function }%
and the function given by $E_{L}=V(L)-L$ as \textit{the para-energy function}
associated with $L$.

The operator $i_{J}$ induced by $J$ and shown as
\begin{equation}
i_{J}\omega (Z_{1},Z_{2},...,Z_{r})=\sum_{i=1}^{r}\omega
(Z_{1},...,JZ_{i},...,Z_{r})  \label{2.8}
\end{equation}
is said to be \textit{vertical derivation, }where $\omega \in \wedge ^{r}TM,$
$Z_{i}\in \chi (TM).$ The \textit{vertical differentiation} $d_{J}$ is
defined as follows:
\begin{equation}
d_{J}=[i_{J},d]=i_{J}d-di_{J},  \label{2.9}
\end{equation}
where $d$ is the usual exterior derivation. For almost paracomplex structure
$J$ determined by (\ref{2.4}), the closed para-K\"{a}hlerian form is the
closed 2-form given by $\Phi _{L}=-dd_{J}L$ such that
\begin{equation}
d_{J}=-\mathbf{j}\frac{\partial }{\partial z^{i}}dz^{i}+\mathbf{j}\frac{%
\partial }{\partial \overline{z}^{i}}d\overline{z}^{i}:\mathcal{F}%
(TM)\rightarrow \wedge ^{1}TM.  \label{2.10}
\end{equation}
\textit{Paracomplex-Euler-Lagrange equations} on
para-K\"{a}hlerian manifold $TM$ are shown by
\begin{equation}
\mathbf{j}\frac{\partial }{\partial t}\left( \frac{\partial L}{\partial z^{i}%
}\right) +\frac{\partial L}{\partial z^{i}}=0,\,\,\,\,\,\,\,\,\mathbf{j}%
\frac{\partial }{\partial t}\left( \frac{\partial L}{\partial \overline{z}%
^{i}}\right) -\frac{\partial L}{\partial \overline{z}^{i}}=0.  \label{2.11}
\end{equation}
Thus, the triple $(TM,\Phi _{L},\xi )$ is called a
\textit{paracomplex-mechanical system}.

\subsection{Paracomplex Hamiltonian Systems}

Here, we consider paracomplex-Hamiltonian equations for classical
mechanics structured on para-K\"{a}hlerian manifold $T^{*}M$
introduced in \cite
{tekkoyun}. Let $T^{*}M$ be any para-K\"{a}hlerian manifold and $(z_{i},%
\overline{z}_{i})$ its coordinates. $\{\frac{\partial }{\partial z_{i}}|_{p},%
\frac{\partial }{\partial \overline{z}_{i}}|_{p}\}$ and $\{\left.
dz_{i}\right| _{p},\left( d\overline{z}_{i}\right) _{p}\}$ be bases over
paracomplex number $\mathbf{A}$ of tangent space $T_{p}(TM)$ and cotangent
space $T_{p}^{*}(TM)$ of $TM.$ Assume that $J^{*}$ is an almost paracomplex
structure given by $J^{*}(dz_{i})=-$\textbf{j}$dz_{i}$, $J^{*}(d\overline{z}%
_{i})=$\textbf{j}$d\overline{z}_{i}$ and $\lambda $ is a para-Liouville form
given by $\lambda =J^{*}(\omega )=\frac{1}{2}$\textbf{j}$(z_{i}d\overline{z}%
_{i}-\overline{z}_{i}dz_{i})$ such that paracomplex 1-form $\omega =\frac{1}{%
2}(z_{i}d\overline{z}_{i}+\overline{z}_{i}dz_{i})$ on $T^{*}M.$ If $\Phi
=-d\lambda $ is closed para-K\"{a}hlerian form$,$ then $\Phi $ is also a
para-symplectic structure on $T^{*}M$.

Let $T^{*}M$\textbf{\ }be para-K\"{a}hlerian manifold with\textbf{\ }closed
para-K\"{a}hlerian form $\Phi .$ Then para-Hamiltonian vector field $Z_{H}$%
\textbf{\ }on\ $T^{*}M$ with\textbf{\ }closed\textbf{\ }form $\Phi $ can be
given by:
\begin{equation}
Z_{H}=-\mathbf{j}\frac{\partial H}{\partial \overline{z}_{i}}\frac{\partial
}{\partial z_{i}}+\mathbf{j}\frac{\partial H}{\partial z_{i}}\frac{\partial
}{\partial \overline{z}_{i}}.  \label{2.12}
\end{equation}
According to (\ref{1.5}), \textit{para}-\textit{Hamiltonian equations} on
para-K\"{a}hlerian manifold $T^{*}M$ are denoted by equations of
\begin{equation}
\frac{dz_{i}}{dt}=-\mathbf{j}\frac{\partial H}{\partial \overline{z}_{i}}%
,\,\,\,\,\,\,\,\,\,\,\,\frac{d\overline{z}_{i}}{dt}=\mathbf{j}\frac{\partial
H}{\partial z_{i}}.  \label{2.13}
\end{equation}

\textbf{Example: }A central force field $f(\rho )=A\rho ^{\alpha
-1}(\alpha \neq 0,1)$ acts on a body with mass $m$ in a constant
gravitational field. Then let us find out the para-Lagrangian and
para-Hamiltonian equations of the motion by assuming the body
always on the vertical plane.

The para-Lagrangian and para-Hamiltonian functions of the system are,
respectively,

\[
L=\frac{1}{2}m\stackrel{.}{z}\stackrel{.}{\overline{z}}-\frac{A}{\alpha }(%
\sqrt{z\overline{z}})^{\alpha }-\mathbf{j}mg\frac{(z-\overline{z})\sqrt{z%
\overline{z}}}{(z+\overline{z})\sqrt{1-\frac{(z-\overline{z})^{2}}{(z+%
\overline{z})^{2}}}},
\]

\[
H=\frac{1}{2}m\stackrel{.}{z}\stackrel{.}{\overline{z}}+\frac{A}{\alpha }(%
\sqrt{z\overline{z}})^{\alpha }+\mathbf{j}mg\frac{(z-\overline{z})\sqrt{z%
\overline{z}}}{(z+\overline{z})\sqrt{1-\frac{(z-\overline{z})^{2}}{(z+%
\overline{z})^{2}}}}.
\]

Then, using (\ref{2.11}) and (\ref{2.13}), the so-called para-Lagrangian and
para-Hamiltonian equations of the motion on the para-mechanical systems, can
be obtained, respectively, as follows:

\[
L1:\,\,\,\,\,\mathbf{j}\frac{\partial }{\partial t}S-S=0,\,\,\,\,\,\,\,\,\,%
\,\,L2:\,\,\,\,\,\,\,\mathbf{j}\frac{\partial }{\partial t}U+U=0,
\]

such that

\begin{eqnarray*}
S &=&-\frac{A}{2z}(\sqrt{z\overline{z}})^{\alpha }-\mathbf{j}\frac{mg(z-%
\overline{z})\overline{z}}{2\sqrt{z\overline{z}}(z+\overline{z})W}-\mathbf{j}%
\frac{mg\sqrt{z\overline{z}}}{(z+\overline{z})W} \\
&&+\mathbf{j}\frac{mg\sqrt{z\overline{z}}(z-\overline{z})}{(z+\overline{z}%
)^{2}W}+\mathbf{j}\frac{mg\sqrt{z\overline{z}}(z-\overline{z})(-\frac{(z-%
\overline{z})}{(z+\overline{z})^{2}}+\frac{(z-\overline{z})^{2}}{(z+%
\overline{z})^{3}})}{(z+\overline{z})W^{3}},
\end{eqnarray*}

\begin{eqnarray*}
U &=&-\frac{A}{2\overline{z}}(\sqrt{z\overline{z}})^{\alpha }-\mathbf{j}%
\frac{mg(z-\overline{z})z}{2\sqrt{z\overline{z}}(z+\overline{z})W}+\mathbf{j}%
\frac{mg\sqrt{z\overline{z}}}{(z+\overline{z})W} \\
&&+\mathbf{j}\frac{mg\sqrt{z\overline{z}}(z-\overline{z})}{(z+\overline{z}%
)^{2}W}+\mathbf{j}\frac{mg\sqrt{z\overline{z}}(z-\overline{z})(\frac{(z-%
\overline{z})}{(z+\overline{z})^{2}}+\frac{(z-\overline{z})^{2}}{(z+%
\overline{z})^{3}})}{(z+\overline{z})W^{3}}
\end{eqnarray*}

and

\begin{eqnarray*}
H1 &:&\,\,\,\,\,\,\,\frac{dz}{dt}=-\mathbf{j(}\frac{A}{2\overline{z}}(\sqrt{z%
\overline{z}})^{\alpha }+\mathbf{j}\frac{mg(z-\overline{z})z}{2\sqrt{z%
\overline{z}}(z+\overline{z})W}-\mathbf{j}\frac{mg\sqrt{z\overline{z}}}{(z+%
\overline{z})W} \\
&&-\mathbf{j}\frac{mg\sqrt{z\overline{z}}(z-\overline{z})}{(z+\overline{z}%
)^{2}W}-\mathbf{j}\frac{mg\sqrt{z\overline{z}}(z-\overline{z})(\frac{(z-%
\overline{z})}{(z+\overline{z})^{2}}+\frac{(z-\overline{z})^{2}}{(z+%
\overline{z})^{3}})}{(z+\overline{z})W^{3}}), \\
H2 &:&\,\,\,\,\,\,\,\,\,\,\,\frac{d\overline{z}}{dt}=\mathbf{j(}\frac{A}{2z}(%
\sqrt{z\overline{z}})^{\alpha }+\mathbf{j}\frac{mg(z-\overline{z})\overline{z%
}}{2\sqrt{z\overline{z}}(z+\overline{z})W}+\mathbf{j}\frac{mg\sqrt{z%
\overline{z}}}{(z+\overline{z})W} \\
&&-\mathbf{j}\frac{mg\sqrt{z\overline{z}}(z-\overline{z})}{(z+\overline{z}%
)^{2}W}-\mathbf{j}\frac{mg\sqrt{z\overline{z}}(z-\overline{z})(-\frac{(z-%
\overline{z})}{(z+\overline{z})^{2}}+\frac{(z-\overline{z})^{2}}{(z+%
\overline{z})^{3}})}{(z+\overline{z})W^{3}}).
\end{eqnarray*}

where $W=\sqrt{1-\frac{(z-\overline{z})^{2}}{(z+\overline{z})^{2}}}.$

\section{Constrained Paracomplex Lagrangians}

In this section, we obtain para-Euler-Lagrange equations with constraints
for classical mechanics structured on para-K\"{a}hlerian manifold $TM$.

Let $J$ be an almost paracomplex structure on the para-K\"{a}hlerian
manifold and $(z^{i},\overline{z}^{i})$ its coordinates. Let us take a
second order differential equation to the vector field $\xi $\ given by:

\begin{equation}
\xi =\xi _{L}+\wedge ^{a}\omega _{a}=\xi ^{i}\frac{\partial }{\partial z^{i}}%
+\overline{\xi }^{i}\frac{\partial }{\partial \overline{z}^{i}}+\wedge
^{a}\omega _{a},\,1\leq a\leq r,  \label{3.1}
\end{equation}
The vector field $V=J\xi _{L}$ calculated by
\begin{equation}
-\mathbf{j}\xi ^{i}\frac{\partial }{\partial z^{i}}+\mathbf{j}\overline{\xi }%
^{i}\frac{\partial }{\partial \overline{z}^{i}},  \label{3.2}
\end{equation}
is \textit{para}-\textit{Liouville vector field }on the para-K\"{a}hlerian
manifold $TM$. The closed 2-form expressed by $\Phi _{L}=-dd_{J}L$ is found
to be:
\begin{eqnarray}
\Phi _{L} &=&-\mathbf{j}\frac{\partial ^{2}L}{\partial z^{j}\partial z^{i}}%
dz^{j}\wedge dz^{i}+\mathbf{j}\frac{\partial ^{2}L}{\partial \overline{z}%
^{j}\partial z^{i}}d\overline{z}^{j}\wedge dz^{i}  \label{3.5} \\
&&-\mathbf{j}\frac{\partial ^{2}L}{\partial z^{j}\partial \overline{z}^{i}}%
dz^{j}\wedge d\overline{z}^{i}-\mathbf{j}\frac{\partial ^{2}L}{\partial
\overline{z}^{j}\partial \overline{z}^{i}}d\overline{z}^{j}\wedge d\overline{%
z}^{i},  \nonumber
\end{eqnarray}
where
\begin{equation}
d_{J}=-\mathbf{j}\frac{\partial }{\partial z^{i}}dz^{i}+\mathbf{j}\frac{%
\partial }{\partial \overline{z}^{i}}d\overline{z}^{i}:\mathcal{F}%
(TM)\rightarrow \wedge ^{1}TM.  \label{3.6}
\end{equation}

If $\xi $ is a second order differential equation defined by (\ref{1.3}),
then we have
\begin{eqnarray}
i_{\xi }\Phi _{L} &=&-\mathbf{j}\xi ^{i}\frac{\partial ^{2}L}{\partial
z^{j}\partial z^{i}}\delta _{i}^{j}dz^{i}+\mathbf{j}\xi ^{i}\frac{\partial
^{2}L}{\partial z^{j}\partial z^{i}}dz^{j}+\mathbf{j}\overline{\xi }^{i}%
\frac{\partial ^{2}L}{\partial \overline{z}^{j}\partial z^{i}}\delta
_{i}^{j}dz^{i}-\mathbf{j}\xi ^{i}\frac{\partial ^{2}L}{\partial \overline{z}%
^{j}\partial z^{i}}d\overline{z}^{j}  \nonumber \\
&&-\mathbf{j}\xi ^{i}\frac{\partial ^{2}L}{\partial z^{j}\partial \overline{z%
}^{i}}\delta _{i}^{j}d\overline{z}^{i}+\mathbf{j}\overline{\xi }^{i}\frac{%
\partial ^{2}L}{\partial z^{j}\partial \overline{z}^{i}}dz^{j}-\mathbf{j}%
\overline{\xi }^{i}\frac{\partial ^{2}L}{\partial \overline{z}^{j}\partial
\overline{z}^{i}}\delta _{i}^{j}d\overline{z}^{i}+\mathbf{j}\overline{\xi }%
^{i}\frac{\partial ^{2}L}{\partial \overline{z}^{j}\partial \overline{z}^{i}}%
d\overline{z}^{j}.  \label{3.7}
\end{eqnarray}
Since closed para-K\"{a}hlerian form $\Phi _{L}$ on $TM$ is para-symplectic
structure, it is obtained
\begin{equation}
E_{L}=-\mathbf{j}\xi ^{i}\frac{\partial L}{\partial z^{i}}+\mathbf{j}%
\overline{\xi }^{i}\frac{\partial L}{\partial \overline{z}^{i}}-L
\label{3.8}
\end{equation}
and hence
\begin{equation}
\begin{array}{ll}
dE_{L}+\wedge ^{a}\omega _{a}= & -\mathbf{j}\xi ^{i}\frac{\partial ^{2}L}{%
\partial z^{j}\partial z^{i}}dz^{j}+\mathbf{j}\overline{\xi }^{i}\frac{%
\partial ^{2}L}{\partial z^{j}\partial \overline{z}^{i}}dz^{j}-\frac{%
\partial L}{\partial z^{j}}dz^{j} \\
& -\mathbf{j}\xi ^{i}\frac{\partial ^{2}L}{\partial \overline{z}^{j}\partial
z^{i}}d\overline{z}^{j}+\mathbf{j}\overline{\xi }^{i}\frac{\partial ^{2}L}{%
\partial \overline{z}^{j}\partial \overline{z}^{i}}d\overline{z}^{j}-\frac{%
\partial L}{\partial \overline{z}^{j}}d\overline{z}^{j}+\wedge ^{a}\omega
_{a}.
\end{array}
\label{3.9}
\end{equation}
According to (\ref{1.3}), if (\ref{3.7}) and (\ref{3.9}) are equal to each
other, then the following equation can be obtained:
\begin{equation}
\begin{array}{l}
+\mathbf{j}\xi ^{i}\frac{\partial ^{2}L}{\partial z^{j}\partial z^{i}}dz^{j}+%
\mathbf{j}\overline{\xi }^{i}\frac{\partial ^{2}L}{\partial \overline{z}%
^{j}\partial z^{i}}dz^{j}+\frac{\partial L}{\partial z^{j}}dz^{j} \\
-\mathbf{j}\xi ^{i}\frac{\partial ^{2}L}{\partial z^{j}\partial \overline{z}%
^{i}}d\overline{z}^{j}-\mathbf{j}\overline{\xi }^{i}\frac{\partial ^{2}L}{%
\partial \overline{z}^{j}\partial \overline{z}^{i}}d\overline{z}^{j}+\frac{%
\partial L}{\partial \overline{z}^{j}}d\overline{z}^{j}=\wedge ^{a}\omega
_{a}
\end{array}
\label{3.10}
\end{equation}
Now, let curve $\alpha :\mathbf{A}\rightarrow TM$ be integral curve of $\xi
, $ which satisfies equations of

\begin{equation}
\begin{array}{l}
+\mathbf{j}\left[ \xi ^{j}\frac{\partial ^{2}L}{\partial z^{j}\partial z^{i}}%
+\stackrel{.}{\xi }^{i}\frac{\partial ^{2}L}{\partial \stackrel{.}{z}%
^{j}\partial z^{i}}\right] dz^{j}+\frac{\partial L}{\partial z^{j}}dz^{j} \\
\\
-\mathbf{j}\left[ \xi ^{j}\frac{\partial ^{2}L}{\partial z^{j}\partial
\stackrel{.}{z}^{i}}+\stackrel{.}{\xi }^{j}\frac{\partial ^{2}L}{\partial
\stackrel{.}{z}^{j}\partial \stackrel{.}{z}^{i}}\right] d\stackrel{.}{z}^{j}+%
\frac{\partial L}{\partial \stackrel{.}{z}^{j}}d\stackrel{.}{z}^{j}=\wedge
^{a}\omega _{a},
\end{array}
\label{3.11}
\end{equation}
where the dots mean derivatives with respect to time and $\omega
_{a}=(\omega _{a})_{i}$ $dz^{i}+(\stackrel{.}{\omega }_{a})_{i}$ $d\stackrel{%
.}{z}^{i}.$

This refers to equations of
\begin{equation}
\frac{\partial L}{\partial z^{i}}+\mathbf{j}\frac{\partial }{\partial t}%
\left( \frac{\partial L}{\partial z^{i}}\right) =\wedge ^{a}(\omega
_{a})_{i},\,\,\,\,\,\,\,\,\,\frac{\partial L}{\partial \stackrel{.}{z}^{i}}-%
\mathbf{j}\frac{\partial }{\partial t}\left( \frac{\partial L}{\partial
\stackrel{.}{z}^{i}}\right) =\wedge ^{a}(\stackrel{.}{\omega }_{a})_{i}.
\label{3.12}
\end{equation}

Thus, the equations obtained in (\ref{3.12}) on para-K\"{a}hlerian manifold $%
TM$ are so-called as \textit{constrained paracomplex
Euler-Lagrange equations}. Then the quartet $(TM,\Phi _{L},\xi
,\overline{\omega })$ is named \textit{constrained paracomplex
mechanical system}.

\section{Constrained Paracomplex Hamiltonians}

Here, we conclude paracomplex Hamiltonian equations with constraints on
para-K\"{a}hlerian manifold $T^{*}M$. Similar to\textbf{\ }(\ref{1.5}),%
\textbf{\ }the vector fields on\textbf{\ }$T^{*}M$ satisfying the condition

\begin{equation}
i_{Za}\Phi =\omega _{a},\,1\leq a\leq s,  \label{4.1}
\end{equation}
can be represented by $Z_{a}.$

\textbf{Proposition: }Let $T^{*}M$\textbf{\ }be para-Kaehlerian manifold with%
\textbf{\ }closed para-K\"{a}hlerian form $\Phi .$ Let us consider vector
field $Z_{a}$\textbf{\ }on\textbf{\ }$T^{*}M$ given by:
\begin{equation}
Z_{a}=-\mathbf{j}(B_{a})_{i}\frac{\partial }{\partial z_{i}}+\mathbf{j}%
(A_{a})_{i}\frac{\partial }{\partial \overline{z}_{i}},  \label{4.2}
\end{equation}

\textbf{Proof:} Let $T^{*}M$\textbf{\ }be para-K\"{a}hlerian manifold with%
\textbf{\ }form $\Phi .$ Consider that vector field $Z_{a}$ is given by
\begin{equation}
Z_{a}=(Z_{a})_{i}\frac{\partial }{\partial z_{i}}+(\overline{Z}_{a})_{i}%
\frac{\partial }{\partial \overline{z}_{i}}.  \label{4.3}
\end{equation}

From (\ref{4.1}), $i_{Z_{a}}\Phi $ can be calculated as
\begin{equation}
i_{Z_{a}}(-d\lambda )=\mathbf{j}(\bar{Z}_{a})_{i}dz_{i}-\mathbf{j}%
(Z_{a})_{i}d\overline{z}_{i}.  \label{4.4}
\end{equation}
Moreover, we set
\begin{equation}
\omega _{a}=(A_{a})_{i}dz_{i}+(B_{a})_{i}d\overline{z}_{i}  \label{4.5}
\end{equation}
According to (\ref{4.1})$,$ if (\ref{4.4}) and (\ref{4.5}) are equal to each
other, proof finishes$.$ $\mathbf{\diamondsuit }$

Now, with the case of (\ref{1.5}) and (\ref{1.7}) and (\ref{4.1}); one may
easily deduce

\begin{equation}
Z=Z_{H}+\wedge ^{a}Z_{a}.  \label{4.6}
\end{equation}
Hence, by means of (\ref{2.13}), (\ref{4.2}) and (\ref{4.6}) we obtain the
following vector field
\begin{equation}
Z=-\mathbf{j}(\frac{\partial H}{\partial \overline{z}_{i}}+\wedge
^{a}(B_{a})_{i})\frac{\partial }{\partial z_{i}}+\mathbf{j}(\frac{\partial H%
}{\partial z_{i}}+\wedge ^{a}(A_{a})_{i})\frac{\partial }{\partial \overline{%
z}_{i}}.  \label{4.7}
\end{equation}
Suppose that curve
\[
\alpha :I\subset \mathbf{A}\rightarrow T^{*}M
\]
be an integral curve of paracomplex vector field $Z$ given by (\ref{4.7})$,$
i.e.,
\begin{equation}
Z(\alpha (t))=\stackrel{.}{\alpha }(t),\,\,t\in I.  \label{4.8}
\end{equation}
In the local coordinates, for $\alpha (t)=(z_{i}(t),\overline{z}_{i}(t)),$
we have

\begin{equation}
\stackrel{.}{\alpha }(t)=\frac{dz_{i}}{dt}\frac{\partial }{\partial z_{i}}+%
\frac{d\overline{z}_{i}}{dt}\frac{\partial }{\partial \overline{z}_{i}}.
\label{4.9}
\end{equation}
Then we reach the following equations

\begin{equation}
\begin{array}{l}
\frac{dz_{i}}{dt}=-\mathbf{j}(\frac{\partial H}{\partial \overline{z}_{i}}%
+\wedge ^{a}(B_{a})_{i}), \\
\frac{d\overline{z}_{i}}{dt}=\mathbf{j}(\frac{\partial H}{\partial z_{i}}%
+\wedge ^{a}(A_{a})_{i}), \\
(A_{a})_{i}\frac{dz_{i}}{dt}+(B_{a})_{i}\frac{d\overline{z}_{i}}{dt}=0,
\end{array}
\label{4.10}
\end{equation}
which are so-called as \textit{constrained paracomplex Hamiltonian
equations} on para-K\"{a}hlerian manifold $T^{*}M$. Here $1\leq
a\leq s.$ Then the quartet $(T^{*}M,\Phi,H ,\overline{\omega })$
is named \textit{constrained paracomplex mechanical system}.

\section{Conclusion}

Finally, considering the above, complex analogous of the geometrical and
mechanical meaning of constraints given in \cite
{weber,deleon,tekkoyun1,kiehn} can be explained as follows:

\textbf{1) }Let $\overline{\omega }$ be a system of constraints on
para-K\"{a}hlerian manifold $TM$ or $T^{*}M.$ Then it may be defined a
distribution $D$ or $D^{*}$on $\overline{\omega }$ as follows:
\begin{equation}
\begin{array}{l}
D(x)=\{\left. \xi \in T_{x}TM\right| \,\omega _{a}(\xi )=0,\,\mbox{ for all }%
\;a,\,1\leq a\leq r\} \\
D^{*}(x)=\{\left. Z\in T_{x}T^{*}M)\right| \,\omega _{a}(Z)=0,\,%
\mbox{ for
all }\;a,1\leq a\leq s\}
\end{array}
\label{5.1}
\end{equation}
Thus $D$ or $D^{*}$ is $(2m-r)$ or $(2m-s)$-dimensional distribution on $TM$
or $T^{*}M.$ In this case, a system of paracomplex constraints $\overline{%
\omega }$ is paraholonomic, if the distribution $D$ or $D^{*}$ is
integrable; otherwise $\overline{\omega }$ is paraanholonomic. \thinspace
Hence, $\overline{\omega }$ is paraholonomic if and only if the ideal $\rho $
of $\wedge TM$ or $\wedge T^{*}M$ generated by $\overline{\omega }$ is a
differential ideal, i.e., $d\rho \subset \rho .$ Obviously, (\ref{3.12}) and
(\ref{4.10}) hold both paraholonomic and paraanholonomic constraints. The
motion for a system of paraholonomic constraints lies on a specific leaf of
the foliation defined by $D$ or $D^{*}.$

\textbf{2) }From (\ref{1.3}) and (\ref{1.7}), the following equalities can
be obtained:
\begin{equation}
\begin{array}{l}
0=(i_{\xi }\Phi )(\xi )=dE_{L}(\xi )=\xi (E_{L}), \\
0=(i_{Z}\omega )(Z)=dH(Z)=Z(H).
\end{array}
\label{5.2}
\end{equation}
So, Lagrangian energy $E_{L}$\ and Hamiltonian energy $H$ of (\ref{3.12})
and (\ref{4.10}) for a solution $\alpha (t)$ are, respectively, conserved.

\textbf{Acknowledgment}

We are very grateful to Professor M. Adak for his valuable comments, remarks
and suggestions.

\end{titlepage}
\end{document}